\documentclass[aps,pre,preprint,onecolumn,groupedaddress,superscriptaddress,showpacs,nofootinbib,notitlepage]{revtex4-1}

\usepackage{hyperref}       
\usepackage{url}            
\usepackage{booktabs}       
\usepackage{amsfonts}       
\usepackage{amssymb}
\usepackage{mathtools}
\usepackage{nicefrac}       
\usepackage{microtype}      
\usepackage{amsmath}
\usepackage{xspace}
\usepackage{xcolor,colortbl}
\usepackage{hhline}
\usepackage{bbm}
\usepackage{changes}
\usepackage{tabularx}
\usepackage{amsthm}
\usepackage{multirow, makecell}
\usepackage{slashbox}
\usepackage{changes}
\usepackage{tikz}
\usepackage{color,soul}
\usepackage{framed} 
\usepackage{multirow}
\usepackage{comment}

\newtheorem{theorem}{Theorem}

\begin{document}

\title{Tradeoffs in Hierarchical Voting Systems}
\author{Lucas B\"ottcher}
\email{l.boettcher@fs.de}
\affiliation{Computational Social Science, Frankfurt School of Finance and Management, Frankfurt am Main, 60322, Germany}
\affiliation{Dept.~of Computational Medicine, University of California, Los Angeles, 90095-1766, Los Angeles, CA, United States}
\author{Georgia Kernell}
\email{gkernell@ucla.edu}
\affiliation{Depts.~of Communication and Political Science, University of California, Los Angeles, 90095-1766, Los Angeles, CA, United States}
\date{\today}
\begin{abstract}
Condorcet's jury theorem states that the correct outcome is reached in direct majority voting systems with sufficiently large electorates as long as each voter's independent probability of voting for that outcome is greater than 0.5. Yet, in situations where direct voting systems are infeasible, such as due to high implementation and infrastructure costs, hierarchical voting systems provide a reasonable alternative. We study differences in outcome precision between hierarchical and direct voting systems for varying group sizes, abstention rates, and voter competencies. Using asymptotic expansions of the derivative of the reliability function (or Banzhaf number), we first prove that indirect systems differ most from their direct counterparts when group size and number are equal to each other, and therefore to $\sqrt{N_{\rm d}}$, where $N_{\rm d}$ is the total number of voters in the direct system. In multitier systems, we prove that this difference is maximized when group size equals $\sqrt[n]{N_{\rm d}}$, where $n$ is the number of hierarchical levels. Second, we show that while direct majority rule always outperforms hierarchical voting for homogeneous electorates that vote with certainty, as group numbers and size increase, hierarchical majority voting gains in its ability to represent all eligible voters. Furthermore, when voter abstention and competency are correlated within groups, hierarchical systems often outperform direct voting, which we show by using a generating function approach that is able to analytically characterize heterogeneous voting systems.
\end{abstract}
\maketitle
\pagebreak
\section*{Introduction}
Condorcet’s jury theorem~\cite{de2014essai} characterizes the outcome in information aggregation problems arising in a group of $N$ decision makers (e.g., voters, jurists, or shareholders) that choose by majority vote between a correct outcome with probability $\epsilon$ and an incorrect outcome with probability $1-\epsilon$. The theorem demonstrates that when $\epsilon>1/2$, the correct outcome is reached as $N\rightarrow \infty$. Alternatively, if $\epsilon<1/2$, the chance of achieving the correct outcome is decreasing in $N$. In this case, the probability of obtaining the correct outcome is maximized when there is only one decision maker. When $\epsilon=1/2$ the probability of a correct outcome is 1/2 for any $N$. Note that the outlined theorem is a direct consequence of the weak law of large numbers~\cite{kalai2006threshold}. 

Several remarks are in order to clarify the simplifying assumptions made in Condorcet's formulation of information aggregation. First, Condorcet's theorem assumes a homogeneous pool of voters---both with respect to preferences and accuracy. Everyone votes for the (mutually agreed upon) ``correct'' outcome with equal probability $\epsilon$. Second, it neglects the dependence between voters' decisions and dynamical social-influence effects~\cite{bottcher2017critical}. Third, in Condorcet's setting, majority rule is applied to all decision makers at once; the decision-making process is direct and not based on the voting outcome in different subgroups. Fourth, the model assumes no abstention (or assumes abstention is unrelated to $\epsilon$). Everyone votes with equal probability.

The influence of certain aspects of heterogeneity, voter dependence, and hierarchical voting are discussed in \cite{boland1989majority}. Heterogeneity in information can be captured by substituting the mean voter's opinion for $\epsilon$~\cite{boland1989majority}. In addition, positive correlations between voters' decisions have a negative impact on the effectiveness of direct majority voting systems~\cite{boland1989majority}, a finding confirmed more recently in \cite{kaniovski2011optimal}. Previous work has also shown that the probability of accepting the correct decision under Condorcet's assumptions is larger in direct voting systems than in certain indirect systems~\cite{boland1989majority,boland1989modelling}. As discussed in \cite{berg1998collective}, indirect voting systems may outweigh direct systems when there is a tradeoff between effectiveness (or outcome precision) and implementation costs, which may increase with both group number and size. \footnote{We use the terms ``hierarchical'' and ``indirect'' interchangeably to refer to systems where preferences are first aggregated within subgroups and then aggregated across groups. A formal definition follows in the description of the model.}

Hierarchical majority-voting models are often used to simplify descriptions of representative voting systems~\cite{berg1997indirect}. One marked difference between direct and indirect voting systems is the different minimum number of pivotal voters~\cite{haggstrom2006law,kalai2006threshold}. For example, consider a population of $N=15$ voters, for whom 7 voters prefer party A and 8 voters prefer party B. In a direct voting system, convincing one B voter to switch to A will alter the election outcome. In an indirect voting system with 3 groups of 5 voters each, it is possible that party B wins the election with 6 votes instead of 7 votes if 2 groups each contain 3 Party B supporters. Thus, party A may have to convince 2 additional voters instead of 1.\footnote{Note, however, that as $N\rightarrow \infty$ within groups, the same result is obtained.}

According to May's theorem~\cite{may1952set}, majority voting is the only voting rule that satisfies certain fairness properties. Alternative rules, such as unanimity voting, are often inferior to majority voting because they confer larger type I and type II errors~\cite{feddersen1999elections}. For example, Feddersen and Pesendorfer~\cite{feddersen1998convicting} find that unanimity rule leads juries to be more likely to both convict an innocent person and acquit a guilty person.\footnote{Unanimous voting decisions are also required in various other settings, such as to suspend certain European Union membership rights according to article 7 of the Treaty on European Union~\cite{article7}.} For social choice problems with at least three voting alternatives, Arrow's impossibility theorem~\cite{arrow1950difficulty,miller2019reflections,sen2020majority} states that there exists no social welfare function (or rank-order electoral system) that maps individual preferences to a societal preference order, satisfying desired fairness conditions. As a further generalization of Condorcet's direct binary choice voting scheme, majority runoff elections with three candidates and a continuous range of voter preferences were modeled in \cite{bouton2015majority}. In such multi-stage elections, the Condorcet winner (i.e., the candidate who would receive the majority of votes in a one-to-one contest against any other candidate) may not even participate in the final runoff round~\cite{bouton2015majority}.

Majority rules and hierarchical decision making are relevant not only in social choice problems~\cite{brandt2016handbook} and organizational decision making~\cite{csaszar2013organizational,gersbach2020appointed}, but also to describe the interactions between grid cells in cellular automata~\cite{gartner2017color} and to perform renormalization operations in statistical mechanics~\cite{CUP_StatPhys}. Another area of application of majority functions is information aggregation in neural networks~\cite{richards2006neural} and, in particular, ensemble machine learning, where the outcomes of individual classifiers are combined using (weighted) majority voting to achieve better model performances~\cite{dietterich2000ensemble}. Thus, identifying the factors that impact outcome precision in hierarchical voting systems is crucial for our broad understanding of preference and information aggregation across fields.

In this paper, we examine differences in indirect and direct voting systems, focusing on variations in group size, group number, abstention rates and voter competencies. First, we introduce a simple model of hierarchical voting and show that for any number of voters, groups, and $\epsilon>0.5$, direct voting models outperform indirect models. Simply put, indirect voting introduces room for error. Although under certain conditions indirect systems approach the same outcome as in the direct voting model, adding layers of aggregation is always strictly dominated in expectation by a direct voting model.

Next, we solve for the probability the indirect model leads to the correct outcome as a function of several parameters. We show that hierarchical voting increases in accuracy as (a) the number of levels, (b) the number of groups, and (c) the number of voters within groups increases, all else equal, for $\epsilon>0.5$, where $\epsilon$ equals the probability a given individual votes for the correct outcome. (The opposite is true for $\epsilon<0.5$.) Earlier results on hierarchical voting systems with two layers~\cite{boland1989majority} led to the conjecture that the collective performance of a group of independent voters is greater for a large number of small groups than for a small number of large groups. A more recent study~\cite{kao2019modular} provided numerical evidence that the difference between a two-tier indirect and direct system is greatest when the number of groups is similar to the number of members per group. Using an asymptotic expansion of the derivative of the reliability function (or Banzhaf number~\cite{berg1997indirect}), we solve this long-standing debate by proving that the outcome for any indirect system with $n$ tiers and $N_{\rm d}$ voters differs most from a direct voting system when group size and number equal $\sqrt[n]{N_{\rm d}}$.

The next section studies heterogeneous voting systems and how voter abstention may impact (and even change) the superiority of direct versus indirect voting. To describe heterogeneous multi-tier voting systems analytically, we develop a corresponding generating function approach. To the best of our knowledge, this is the first study that uses such a formalism to model systems with heterogeneous voter competencies. Previous work has relied on different bounds~\cite{hoeffding1956distribution,hodges1960poisson,boland1989majority} to characterize such systems, while the proposed generating function approach has the benefit that it directly captures underlying heterogeneities without relying on further approximations. We find that when the likelihood of voting is uncertain but equal across a population, direct voting remains superior. However, when abstention and voter competence are heterogeneous across subgroups -- as is likely the case in many electoral and managerial contexts -- indirect voting gains in its ability to represent the entire electorate. Indeed, in some cases, hierarchical voting leads to dramatically different outcomes. Another way to explain this finding is that when a voting system is meant to represent the preferences of \emph{eligible}, rather than actual, voters, indirect elections  provide a potential correction to direct voting.
\section*{A simple model of hierarchical voting}
\begin{figure*}
    \centering
    \includegraphics[width = \textwidth]{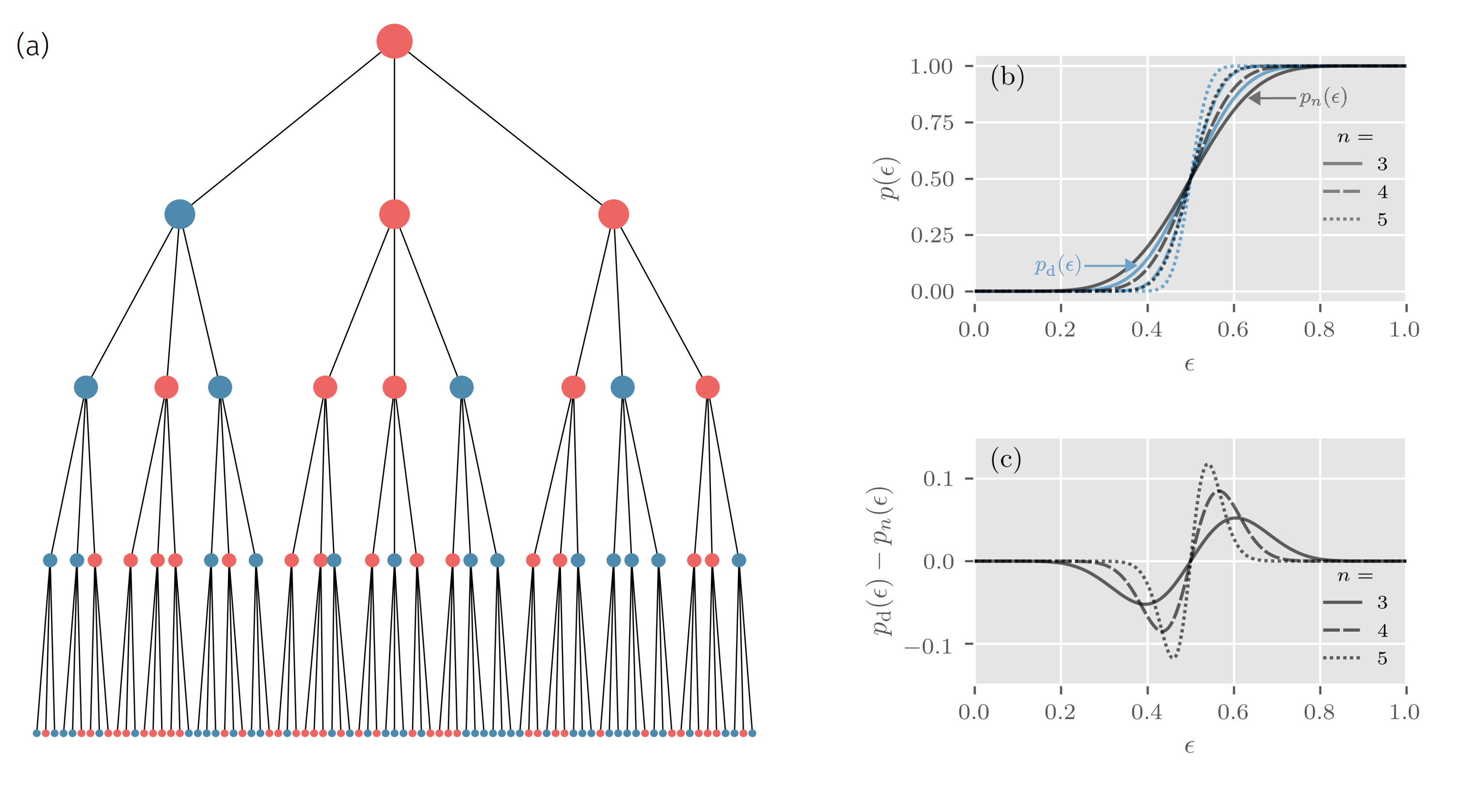}
    \caption{\textbf{Hierarchical and direct majority voting.} (a) Example of a hierarchical voting system with four layers. Blue and red nodes are in states $0$ and $1$, respectively. (b) Proportion of correct voting outcomes in three direct voting systems (denoted by $p_d(\epsilon)$ in the text) and three indirect voting systems (denoted by $p_n(\epsilon)$ in the text) as a function of the probability of voting for the correct outcome (denoted by $\epsilon$ in the text). Different curves correspond to different numbers of layers $n$. In the shown examples, we set $n=3,4,5$ and the number of voters in each majority-voting group is $k=3$. (c) Difference between $p_{\rm d}(\epsilon)$, the proportion of correct voting outcomes in direct voting systems, and $p_n(\epsilon)$ as a function of $\epsilon$. The number of voters in the direct voting system is $N_{\rm d}=k^n$. All parameters are as in (b). In the vicinity of the Condorcet threshold $\epsilon=0.5$, the slopes of the functions describing the voting outcome of direct majority voting are steeper compared those observed in hierarchical voting systems with the same number of voters in the bottom layer, illustrating that direct majority systems lead to the correct outcome more often if $\epsilon>0.5$.}
    \label{fig:problem}
\end{figure*}
We consider two voting systems: one with indirect, or hierarchical, voting and another with direct voting. The hierarchical voting system we study is represented by a regular tree network $G(n,k)$ with $n$ branches. Each node (except leaf nodes which represent the ``electorate'') has $k$ child nodes, where $k$ is an odd number. The total number of nodes in such an indirect majority-voting  model is $N=\sum_{i=0}^{n} k^i=(k^{n+1}-1)/(k-1)$ (including the top node). In comparison, the number of voters in a direct voting system is $N_{\rm d}=k^n$. Nodes can be in two different states: ``$0$'' and ``$1$'', representing two different policy positions. Initially, we set each node to $1$ with probability $\epsilon$. Next, we apply a recursive majority rule~\cite{haggstrom2006law} to update the states of all root nodes. That is, the state of a node is set to $0$ if the majority of its children are in state $0$. Otherwise, it will be set to $1$. Note that a tie cannot occur because $k$ is an odd number. Figure~\ref{fig:problem}(a) shows an example of a hierarchical voting system with $n=4$ layers.

To characterize a hierarchical voting system mathematically, we denote by $p_n(\epsilon)$ and $1-p_n(\epsilon)$ the probabilities that the top node is in state $1$ and $0$, respectively. According to Condorcet's theorem~\cite{de2014essai}, we find
\begin{equation}
\lim_{n \rightarrow \infty} p_n(\epsilon)=
\begin{cases}
0\,,\quad & \epsilon < 0.5\,,\\
0.5\,,\quad & \epsilon=0.5\,,\\
1\,,\quad &\epsilon > 0.5\,.
\end{cases}
\end{equation}
For systems of finite size, the probability $p_n(\epsilon)$ does not jump from 0 to 1 at $\epsilon=0.5$, but instead smoothly approaches 1 for large enough values of $\epsilon$ [Fig.~\ref{fig:problem}(b)]. The probability, $p_n(\epsilon)$, that the top node of such a finite system is in state 1 can be defined recursively in terms of the model parameters. Let $X_i$ be a binomial random variable representing the number of nodes in state 1 in stage $i$ $(1 \leq i \leq n)$. That is, $X_i \sim B(k,p_i)$ where $p_i$ is the probability that a node in stage $i$ is in state 1. For fixed values of $k$, $n$, and $p_0(\epsilon)=\epsilon$, then, the recursion relation is given by
\begin{equation}
p_i(\epsilon) = \Pr \left(X_{i-1} > \dfrac{k}{2}\right)=\sum\limits_{j=\lceil \frac{k}{2} \rceil}^{k} {k \choose j} p^j_{i-1}{(1-p_{i-1})}^{k-j}\,,
\label{eq:recursive}
\end{equation}
where $\lceil k/2 \rceil$ is the ceiling above $k/2$ (i.e., the smallest integer greater than or equal to $k/2$). 

For a comparison with a direct voting system, we denote by $p_{\rm d}(\epsilon)$ and $1-p_{\rm d}(\epsilon)$ the corresponding probabilities of reaching a correct and incorrect voting outcome. Similarly to Eq.~\eqref{eq:recursive}, the reliability function $p_{\rm d}(\epsilon)$ is given by
\begin{equation}
p_{\rm d}(\epsilon) = \Pr \left(X > \dfrac{N_{\rm d}}{2}\right)=\sum\limits_{j=\lceil \frac{N_{\rm d}}{2} \rceil}^{N_{\rm d}} {N_{\rm d} \choose j} \epsilon^j{(1-\epsilon)}^{N_{\rm d}-j}\,,
\end{equation}
where $X$ denotes the number of voters in state $1$. In the limit $N_{\rm d}\rightarrow \infty$ and $N\rightarrow\infty$, it holds that $p_{\rm d}(\epsilon)=p_{n}(\epsilon)$. This equality is approximately fulfilled if the number of voters is sufficiently large. In Fig.~\ref{fig:problem}(c), we show the difference $p_{\rm d}(\epsilon)-p_n(\epsilon)$ for different values of $n$ and observe that the direct voting system outperforms the hierarchical one since $p_{\rm d}(\epsilon)>p_n(\epsilon)$ for $\epsilon > 0.5$ and $p_{\rm d}(\epsilon)<p_n(\epsilon)$ for $\epsilon < 0.5$. In other words, direct voting systems lead to the correct voting outcome more often than hierarchical ones if the probability of voting for the correct outcome is larger than $0.5$. We also observe in Fig.~\ref{fig:problem}(c) that differences between direct and hierarchical voting systems vanish for sufficiently small and large values of $\epsilon$. The greater the number of layers $n$, the more pronounced the differences in the vicinity of $\epsilon\approx 0.5$. 
\begin{figure}
    \centering
    \includegraphics[width=\textwidth]{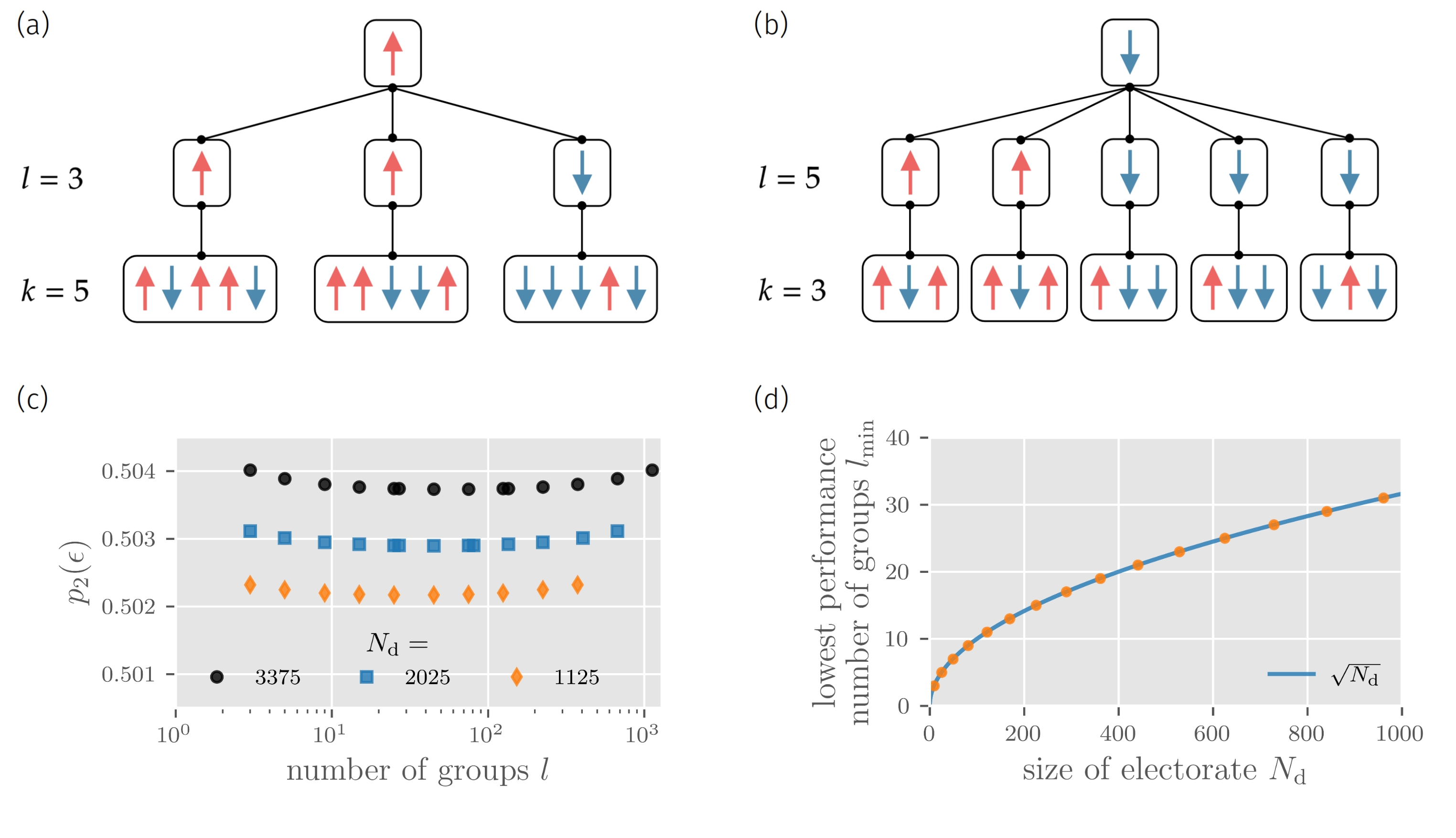}
    \caption{\textbf{Effect of different group sizes on performance of two-tier voting systems.} (a) A two-tier electoral system with $l=3$ groups and $k=5$ voters within each group. Arrows pointing upwards (red) and downwards (blue) represent individuals who vote for the correct and incorrect outcome, respectively. (b) For the same distribution of voter types in the electorate, the voting outcome is different than in (a) because of the different voting group structure. (c) The probability $p_2(\epsilon)$ of a correct outcome for different sized electorates, $N_{\rm d}=kl$, and  $\epsilon=0.5001$ as a function of the number of groups $l$. For two of the selected electorate sizes, $N_{\rm d}=1125$ and $3375$, the corresponding square roots are not integer values and approximately 34 and 58, respectively. Hence, the minimum of $p_2(\epsilon)$ is attained for group sizes $l_{\rm min}=25$ and $45$ that are close but not equal to $\sqrt{N_{\rm d}}$. For $N_{\rm d}=2025$, the minimum of $p_2(\epsilon)$ is attained at $l_{\rm min}=\sqrt{N_{\rm d}}= 45$. (d) Orange dots indicate group sizes with the lowest performance, $l_{\rm min}$, for all electorates with integer square roots. The blue solid line describes the number of groups with the lowest performance according to $l_{\rm min}=\sqrt{N_{\rm d}}$.}
    \label{fig:group_sizes}
\end{figure}
In the following theorem, we will formalize and prove these observations for hierarchical majority-voting systems with $n > 1$ layers.

\begin{theorem}
\label{thm:direct_vs_hierarchical}
Let $p_{\rm d}(\epsilon)$ and $p_n(\epsilon)$ be the probabilities of reaching the correct voting outcome in direct and hierarchical majority-voting systems with the same number of voters in the electorate who vote for the correct voting outcome with probability $\epsilon$. For a hierarchical voting system with at least $n=2$ layers, it holds that
\begin{enumerate}
    \item[(i)] $p_{\rm d}(\epsilon) = p_{n}(\epsilon) \quad \text{for} \quad \epsilon \in \{ 0,0.5,1\}$,
    \item[(ii)] $p_{\rm d}(\epsilon) < p_{n}(\epsilon) \quad \text{for} \quad \epsilon \in (0,0.5)$,
    \item[(iii)] $p_{\rm d}(\epsilon) > p_{n}(\epsilon) \quad \text{for} \quad \epsilon \in (0.5,1)$.
\end{enumerate}
(Points (i--iii) imply that the probability of accepting the correct voting outcome is larger in direct majority systems compared to hierarchical ones as long as $\epsilon \in (0.5,1)$.)
\end{theorem}

The proof of \ref{thm:direct_vs_hierarchical} is in Appendix \ref{app:direct_vs_hierarchical}. For indirect voting systems with one layer, similar results are presented in \cite{boland1989majority,boland1989modelling}, where the authors consider an indirect voting system with $n_1$ voter groups, each of size $n_2$ ($n_1$ and $n_2$ being odd integers). Our proof applies to hierarchical voting systems with $n$ layers, making it more general than the proof presented in \cite{boland1989modelling}. Moreover, we also quantify the slope of the reliability functions $p_{\rm d}(\epsilon)$ and $p_{n}(\epsilon)$ at $\epsilon=0.5$.
\section*{Influence of group size and number of groups}
In a hierarchical voting system, the probability of reaching the correct outcome depends not only on the number of layers but also on the voting group size and the number of groups. In a two-tier voting system with $l$ groups and $k$ voters per group [Fig.~\ref{fig:group_sizes}], we have
\begin{equation}
p_1(\epsilon)=\Pr \left(X_{1} > \dfrac{k}{2}\right)=\sum\limits_{j=\lceil \frac{k}{2} \rceil}^{k} {k \choose j} \epsilon^j{(1-\epsilon)}^{k-j}
\label{eq:p1_2tier}
\end{equation}
and
\begin{equation}
p_2(\epsilon)=\Pr \left(X_{2} > \dfrac{l}{2}\right)=\sum\limits_{j=\lceil \frac{l}{2} \rceil}^{l} {l \choose j} p_1^j (1-p_1)^{l-j}\,
\label{eq:p2_2tier}
\end{equation}
where $p_1(\epsilon)$ and $p_2(\epsilon)$ refer to the probabilities of reaching the correct outcomes in layers $1$ and $2$, respectively. In Fig.~\ref{fig:group_sizes}(a,b), we show two two-tier voting systems with $(k,l)=(5,3)$ and $(k,l)=(3,5)$, respectively. Although the distribution of voters in both electorates is identical, the voting outcome differs because of the different voting group structure. 

The probability $p_2(\epsilon)$ of reaching a correct voting outcome increases with the number of voters in the electorate $N_{\rm d}=kl$ for $\epsilon>0.5$. If the number of voters $N_{\rm d}$ is constant, what is the most effective and least effective composition of voters per group $k$ and numbers of groups $l$? Mathematically, for constant $k,l$, $l>k$, and $\epsilon>0.5$, the probability $p_2(\epsilon)$ is maximized if one uses $l$ groups with $k$ members each~\cite{berg1997indirect}. The opposite holds for $l<k$. Hence, the reliability function $p_2(\epsilon)$ is not symmetric in $k,l$. Earlier related observations \cite{boland1989majority} led to the conjecture that the collective performance of a group of independent voters is larger for a large number of small groups than for a small number of large groups. Indeed, next to direct voting, the most effective two-tier voting system is that with the smallest possible number of voters per group $k$ since it resembles direct voting most closely. However, this observation does not imply that the least effective voting system is that with the maximum possible number of voters per group for constant $N_{\rm d}$. Other work has suggested that the difference between indirect and direct voting models is greatest when group size and number are similar~\cite{kao2019modular}. Kao and Couzin (2019) state that ``the modular structure that leads to the lowest collective accuracy occurs very close to when there are $\sqrt{N_{\rm d}}$ subgroups with $\sqrt{N_{\rm d}}$ individuals per subgroup''~\cite{kao2019modular}.

Using an asymptotic expansion of the derivative of the reliability function about one of its fixed points, Theorem \ref{thm:lowest_performance} states that $p_n(\epsilon)$ is minimized when $k=l$. 
\begin{theorem}
\label{thm:lowest_performance}
For $\epsilon\in(0.5,1)$, the probability of reaching a correct voting outcome in a two-tier voting system is minimized asymptotically if the number of groups approaches the number of voters per group, i.e.\ if
\begin{equation}
(k_{\rm min},l_{\rm min})=(\sqrt{N_{\rm d}},\sqrt{N_{\rm d}})\,.
\end{equation}
\end{theorem}
In Appendix \ref{app:lowest_accuracy}, we present a proof of Theorem \ref{thm:lowest_performance} and provide numerical evidence that the above square-root relation also holds for small electorates. That is, when group size and number are equal, the indirect model deviates most from the direct model. Thus, for $\epsilon>0.5$, arriving at the ``correct'' decision is least likely when $k=l$. However, when $\epsilon<0.5$, an indirect model with $k=l$ is most likely to produce the correct result.\footnote{Interestingly, previous work has found that the optimal number of representatives is proportional to the square root of the population~\cite{auriol2012optimal}.} 

We also prove in Appendix \ref{app:lowest_accuracy} that Theorem \ref{thm:lowest_performance} can be extended to hierarchical voting systems with $n$ layers and group sizes $k^{(1)},k^{(2)},\dots,k^{(n)}$. The accuracy of such a multi-tier voting system is minimized for $k_{\rm min}^{(1)}=k_{\rm min}^{(2)}=\ldots=k_{\rm min}^{(n)}=N_{\rm d}^{1/n}$. For example, in a three tier system with $N_{\rm d}=729$ voters, indirect and direct systems deviate most when there are 9 groups of 9 groups of 9 voters, where $\sqrt[3]{729}=9$. In a four-tier system with 6,561 voters, of the 35 possible voting compositions, the difference between direct and indirect systems is maximized for group sizes $\sqrt[4]{6,561}=9$; in an electorate with 2,313,441 voters, out of 471 possible four-tier compositions, the maximum difference occurs at $\sqrt[4]{2,313,441}=39$.

Figure \ref{fig:group_sizes}(c) shows the reliability function $p_2(\epsilon)$ for different numbers of voters in the electorate, $N_{\rm d}$, as a function of group size $l$. We observe that the minimum of $p_2(\epsilon)$ is obtained for group sizes $l$ that are close to $\sqrt{N_{\rm d}}$. For two of the selected electorate sizes ($N_{\rm d}=1125$ and $3375$), the values of $\sqrt{N_{\rm d}}$ do not correspond to integer values. Hence, the group sizes associated with the lowest collective accuracies, $l_{\rm min}$, are close to but not equal to $\sqrt{N_{\rm d}}$. As shown in Fig.~\ref{fig:group_sizes}(d), the relation $l_{\rm min}=\sqrt{N_{\rm d}}$ holds exactly if the square-root of the electorate size $N_{\rm d}$ is a possible group size (i.e., an odd integer).

The intuition behind the $\sqrt{N_d}$ result follows from the law of large numbers. First, in a direct voting system, as $N_d$ increases, the share of voters needed to pass any particular outcome decreases. (For example, if there are 11 total voters, 6 (54.5\%) are needed to vote for an outcome for it to pass. If there are 21 total voters, 11 (52.3\%) are needed to pass an outcome.) Second, a direct system will always require more votes to win than an indirect system with the same number of voters. Consider an election where $N_{\rm d}=81$. In a direct system, at least 41 voters must support any given outcome for it to pass. However, in a two-tiered hierarchical system with 3 groups of 27 voters, the minimum number of voters is equal to 28 (14 in each of 2 groups). (The same minimum number of voters is obtained when there are 27 groups of 3 voters.) And even fewer voters (25) are needed if there are 9 voters in each of 9 groups (where exactly 5 voters favor a measure in each of 5 groups). As the minimum number of voters required to sway an outcome decreases, an indirect voting system deviates even more from a direct system.

In fact, because the probability a correct outcome is achieved increases with the number of required votes to win an election, the fewest voters necessary to achieve any outcome in a two-tier hierarchical system always occurs when group size equals group number. We prove this, as well as the extension to an $n$-tier hierarchical system as stated in Theorem \ref{thm:fewest_voters}, in Appendix \ref{app:fewest_votrs}.

\begin{theorem}
\label{thm:fewest_voters}
The fewest voters necessary to sway an election occurs when $k_{\rm min}^{(i)}=N_{\rm d}^{1/n}$ ($1\leq i \leq n$), where $k_{\rm min}^{(i)}$ is the number of voters (groups) at level $i$.
\end{theorem}
Thus, the probability that a minority of voters sways the overall outcome toward the incorrect choice when $\epsilon>0.5$, or the correct choice when $\epsilon<0.5$, is maximized when the number of groups in each hierarchy equals $N_{\rm d}^{1/n}$. In the standard Condorcet set-up, where all members in society benefit most from the same ``correct'' outcome, hierarchical models with few groups of many voters or many groups of fewer voters outperform those where group size and number are equal. In contrast, in a competitive environment, where voters hold differing preferences over two outcomes (rather than simply differing levels of accuracy), the outcome preferred by a majority of voters may fail if voters supporting the minority opinion are concentrated in a minimum-winning majority of groups. Intentional redistricting to favor one party over another is well documented throughout U.S. history. Theorem \ref{thm:fewest_voters} suggests this partisan gerrymandering may be easier in cases where the number of groups approaches group size.

Moreover, as we show in the next section, when voter participation is not guaranteed, and abstention is positively correlated with accuracy, an indirect model may be more likely to elicit the correct outcome, even in the standard Condorcet set-up.
\section*{Influence of abstention}
Thus far, we have assumed that all voters participate fully in the election. But in any system, people abstain from voting. In some cases, such as either primary or midterm elections in the U.S., less than half of all eligible voters cast their ballot. In this section, we examine how abstention shapes the tradeoffs between direct and indirect voting.

To model abstention in a multi-stage voting system, let $X_i^{(0)}$, $X_i^{(1)}$, and $X_i^{(2)}$ be the number of nodes in layer $i$ that are in state ``0,'' ``1,'' and ``2,'' respectively. Nodes in state ``2'' represent abstaining voters. In a two-tier voting system, uniform abstention (i.e., an equal probability of abstaining for all voters) can be incorporated in Eqs.~\eqref{eq:p1_2tier} and \eqref{eq:p2_2tier} via a multinomial distribution. That is,
\begin{equation}
p_1(\epsilon)=\Pr \left(X_{1}^{(1)} > \dfrac{k}{2}\right)=\sum\limits_{(x_1,x_2,x_3)\in \mathcal{S}(k)} {k \choose x_1, x_2, x_3} \epsilon^{x_1} \alpha^{x_2} (1-\epsilon-\alpha)^{x_3}
\label{eq:multinomial}
\end{equation}
and
\begin{equation}
p_2(\epsilon)=\Pr \left(X_{2}^{(1)} > \dfrac{l}{2}\right)=\sum\limits_{(x_1,x_2,x_3)\in \mathcal{S}(l)} {l \choose x_1, x_2, x_3} p_1^{x_1} \alpha^{x_2} (1-p_1-\alpha)^{x_3}\,,
\end{equation}
where $\mathcal{S}(k)\coloneqq \{(x_1,x_2,x_3)\vert x_1\in\{\lceil \frac{k}{2} \rceil,\dots,k\} \land x_2\in\{0,\dots,k-x_1\} \land x_3=k-x_1-x_2\}$ and $\alpha$ ($0\leq \alpha \leq 1-\epsilon$) is the abstention probability. The condition $0\leq \alpha \leq 1-\epsilon$ guarantees that the term $1-\epsilon-\alpha$ in Eq.~\eqref{eq:multinomial} is positive.

In the above description of (homogeneous) abstention, direct voting is still superior to indirect voting systems since $\alpha$ only reduces the domain of $\epsilon$ from $[0,1]$ to $[0,1-\alpha]$ without altering other previous results.

In real-world scenarios, abstention is often correlated within voting groups. Residents of rural communities may have longer transportation times to reach their polling location, and thus vote at lower rates in the absence of mail-in ballots. Alternatively, voters may be discouraged from voting in urban areas when there are long lines to vote. Turnout is associated with political interest, information, and education~\cite{sondheimer&green2009}, which are not homogeneously distributed across districts. People in economically disadvantaged districts may have fewer financial, temporal, and informational resources available that make voting accessible. Abstention is also related to an election's competitiveness and the number of races on the ballot -- both of which can vary across districts. Even random shocks---like a power outage, hail storm, or freeway accident---can cause geographical correlations in abstention. In a managerial setting, people in certain divisions of a firm may be less informed or interested in company decisions on which they have a vote. And, while all union members can participate and vote in labor decisions, full-time employees, as well as those receiving higher pay and requiring greater skill, are more likely to participate in union votes~\cite{Kolchin&Hyclak1984}. Thus, it is important to examine cases where abstention is not homogeneously distributed across groups.
\begin{figure}
    \centering
    \includegraphics[width=\textwidth]{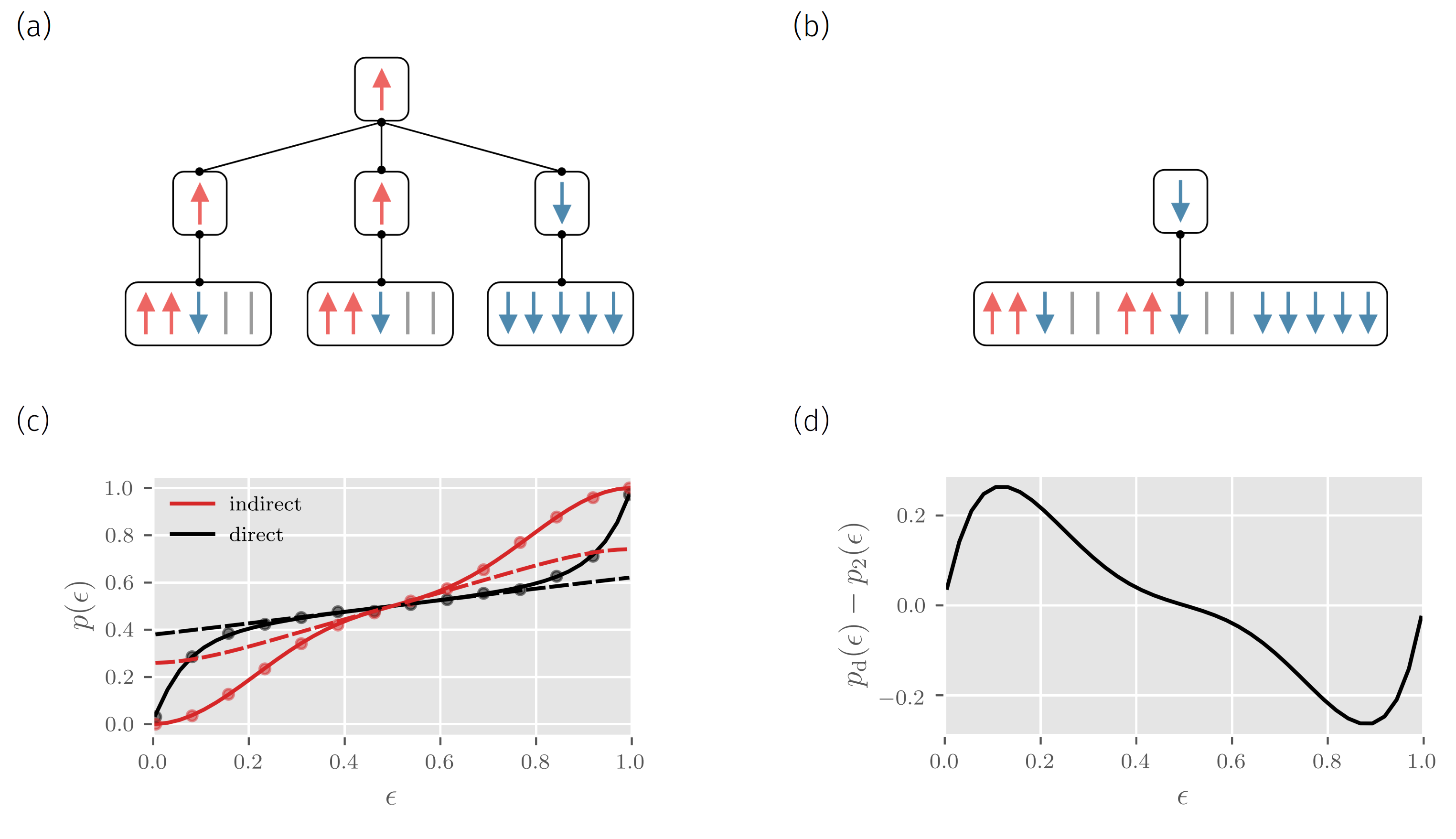}
    \caption{\textbf{Heterogeneous abstention in indirect and direct voting systems.} (a) An indirect voting system with $l=3$ voter groups and voter-group competencies $\epsilon_1=\epsilon_2=\epsilon$ and $\epsilon_3=1-\epsilon$. The number of voters per group are $k_1=k_2=3$ ($\alpha_1=\alpha_2=0.6$) and $k_3=5$ ($\alpha_3=0$), respectively. Solid red and black lines are the analytical solutions Eqs.~\eqref{eq:analytical_abs_1} and \eqref{eq:analytical_abs_2}, respectively. Dashed lines are the corresponding Hoeffding bounds \eqref{eq:abstention_nu_2} and \eqref{eq:hoeffding2}. Red and black disks are numerically obtained voting outcomes that are based on 10,000 samples. (b) The corresponding direct voting system. (c,d) The probability of a correct voting outcome $p(\epsilon)$ for the direct and indirect voting systems as shown in (a,b). The shown data are averages over $3\times10^5$ samples for different values of $\epsilon$. Error bars are smaller than the line width.}
    \label{fig:abstention_nonuniform}
\end{figure}

To model heterogeneous abstention and voter competences within districts, we incorporate the following modifications into our mathematical framework. We again consider a two-tier voting system and account for different abstention probabilities $\alpha_{j}\in[0,1]$ ($1\leq j \leq l$) and voter competences $\epsilon_j$ in each of the $l$ voting groups (or jurisdictions). The number of voters in group $j$ that do not abstain is $\tilde{k}_j=k (1-\alpha_j)$, where $\alpha_j$ is chosen such that $\tilde{k}_j\geq 3$ is an odd number. The probability that the correct voting outcome is achieved in voter group $j$ is
\begin{equation}
p_{1,j}(\epsilon_j)=\Pr \left(X_{1,j}^{(1)} > \dfrac{\tilde{k}_j}{2}\right)=\sum\limits_{m=\lceil \frac{\tilde{k}_j}{2} \rceil}^{\tilde{k}_j} {\tilde{k}_j \choose m} \epsilon_j^m {(1-\epsilon_j})^{\tilde{k}_j-m}\,.
\label{eq:abstention_nu_1}
\end{equation}
According to a theorem by Hoeffding~\cite{hoeffding1956distribution,percus1985probability,boland1989majority}, the probability of a correct voting outcome after information aggregation across all voter groups satisfies
\begin{equation}
p_2(\boldsymbol{\epsilon})=\Pr \left(X_{2}^{(1)} > \dfrac{l}{2}\right)\geq \sum\limits_{j=\lceil\frac{l}{2}\rceil}^l {l \choose j} \bar{p}_{1}^{j} (1-\bar{p}_{1})^{l-j}\,,
\label{eq:abstention_nu_2}
\end{equation}
if the mean voter competence $\bar{p}_1=\bar{p}_1(\boldsymbol{\epsilon})=(p_{1,1}+\dots+p_{1,l})/l \geq 1/2+1/(2 l)$~\cite{boland1989majority} with $\boldsymbol{\epsilon}=(\epsilon_1,\dots,\epsilon_l)$. A similar bound can be given for the probability of a correct voting outcome in direct voting systems:
\begin{equation}
p_{\rm d}(\boldsymbol{\epsilon})=\Pr\left(X>\frac{N_{\rm d}}{2}\right)\geq \sum_{j=\lceil\frac{N_{\rm d}}{2}\rceil}^{N_{\rm d}} {N_{\rm d} \choose j} \bar{\epsilon}^j(1-\bar{\epsilon})^{N_{\rm d}-j}\,,
\label{eq:hoeffding2}
\end{equation}
which holds if $\bar{\epsilon}=(\epsilon_1+\dots+\epsilon_{N_{\rm d}})/N_{\rm d} \geq 1/2+1/(2 N_{\rm d})$. Hoeffding's bound shows that direct heterogeneous voting systems outperform direct homogeneous voting systems with mean voter competency $\bar{\epsilon}$ if the value of $\bar{\epsilon}$ is larger than or equal to $1/2+1/(2 N_{\rm d})$. For the above examples with two levels of hierarchy and three subgroups, the number of voters in the direct voting system is $N_{\rm d}=\tilde{k}_1+\tilde{k}_2+\tilde{k}_3$.

Instead of using the Hoeffding bounds \eqref{eq:abstention_nu_2} and \eqref{eq:hoeffding2} to characterize voting systems with heterogeneous voter competencies and abstention, we use a generating function approach (see Appendix \ref{sec:generating_function}) to obtain an analytical expression for the reliability function of a two-tier voting system. We find that
\begin{equation}
p_2(\boldsymbol{\epsilon})= \frac{1}{P}\left(\frac{p_{1,1} p_{1,2}}{q_{1,1} q_{1,2}}+\frac{p_{1,1} p_{1,3}}{q_{1,1} q_{1,3}}+\frac{p_{1,2} p_{1,3}}{q_{1,2} q_{1,3}}+\frac{p_{1,1} p_{1,2} p_{1,3}}{q_{1,1}q_{1,2}q_{1,3}}\right)\,,
\label{eq:analytical_abs_1}
\end{equation}
where $q_{1,j}=1-p_{1,j}$ and $P=\prod_{j=1}^3 q_{1,j}=q_{1,1}q_{1,2}q_{1,3}$ and that
\begin{equation}
p_{\rm d}(\boldsymbol{\epsilon})= \frac{1}{P}\sum_{j=\lceil\frac{N_{\rm d}}{2}\rceil}^{N_{\rm d}} C_{j}\left(\epsilon_1/(1-\epsilon_1),\dots,\epsilon_{N_{\rm d}}/(1-\epsilon_{N_{\rm d}}) \right)\,,
\label{eq:analytical_abs_2}
\end{equation}
where  $P=\prod_{j=1}^{N_{\rm d}}(1-\epsilon_j)$ and $C_{j}$ is the $j$th elementary symmetric function. As detailed in Appendix \ref{sec:generating_function}, we use Newton's identities to recursively calculate $C_{j}$. Such generating function approaches are useful to analytically capture the properties of general heterogeneous multi-tier voting systems and complement earlier approaches that mainly relied on using different bounds~\cite{hoeffding1956distribution,hodges1960poisson,boland1989majority}.

Contrary to the differences between homogeneous and heterogeneous direct voting systems, as described by the Hoeffding bound \eqref{eq:hoeffding2}, we find that hierarchical \emph{homogeneous} systems outperform \emph{heterogeneous} hierarchical systems. That is, when $\bar{\epsilon}>0.5$, systems where all voters have equal values of $\epsilon$ are associated with higher levels of $p_2$ than those where individuals have heterogeneous values of $\epsilon$, even with a constant mean.

For heterogeneous abstention and homogeneous voter competences (i.e., $\epsilon_j=\epsilon$), direct voting systems are still preferable over indirect ones. However, if high abstention rates are associated with voter groups $j$ in which $\epsilon_j$ is substantially larger than in other groups with low abstention rates, hierarchical voting can outweigh the potential underrepresentation of voter groups with high abstention rates in direct voting systems. Thus, indirect voting can provide an opportunity not only for more efficient preference aggregation, but also for greater representation. In politics, this finding implies greater representation in systems with a high number of representatives (that decide on policy) or delegates (that select a candidate or leader), so long as voter competency is correlated with abstention. If it is not, delegates with fewer representatives should lead to better representation. Parliamentary democracies with plurality or majority legislative elections use indirect systems to select the Prime Minister. All else equal, in countries with high levels of geographical correlation between vote choice and abstention, would expect representation to be greatest when the number of legislators is large (and thus group size is smaller) than when the legislature size is small (and group size is large).



Figure~\ref{fig:abstention_nonuniform}(a) shows an example of a two-tier voting system in which heterogeneous abstention and voter competences impact the final voting outcome. There are $l=3$ voter groups with $k_1=k_2=3$ ($\alpha_1=\alpha_2=0.4$) and $k_3=5$ ($\alpha_3=0$), and the voter-group competencies are $\epsilon_1=\epsilon_2=\epsilon$ and $\epsilon_3=1-\epsilon$, respectively. Figure~\ref{fig:abstention_nonuniform}(b) shows the corresponding direct voting system. For $\epsilon \geq 0.5$, hierarchical voting is associated with a larger probability of a correct voting outcome [Fig.~\ref{fig:abstention_nonuniform}(c)]. For $\epsilon\approx 0.5$, the Hoeffding bounds [dashed lines in Fig.~\ref{fig:abstention_nonuniform}(c)] provide good characterizations of $p_2(\boldsymbol{\epsilon})$ and $p_{\rm d}(\boldsymbol{\epsilon})$. In the limit $\epsilon\rightarrow 1$, the Hoeffding bound of $p_2(\boldsymbol{\epsilon})$ reaches a value of $0.74$ and that of $p_{\rm d}(\boldsymbol{\epsilon})$ reaches a value of about $0.62$, significantly underestimating the true outcome probability. In the limit $\epsilon\rightarrow 1$, both voting systems lead to the correct outcome with probability 1---there are two voters against one voter in the two-tier voting system and 6 against 5 voters in the direct voting system. To verify our analytically obtained voting outcome probabilities \eqref{eq:analytical_abs_1} and \eqref{eq:analytical_abs_2}, we simulated 10,000 independent voting realizations [red and black disks in Fig.~\ref{fig:abstention_nonuniform}(c)] and find excellent agreement between our numerical and analytical results. The maximum difference between the two described voting systems is reached for $\epsilon\approx 0.9$ [Fig.~\ref{fig:abstention_nonuniform}(d)].
\section*{Discussion and conclusion}
This paper investigates differences in direct and indirect voting systems for varying voter competencies, group sizes, and abstention rates. For homogeneous competencies and abstention rates across voter groups, direct voting always outperforms hierarchical voting system. The greater the number of hierarchical layers, the more pronounced the differences between hierarchical and direct voting at the Condorcet threshold $\epsilon=0.5$ (see the proof of Theorem~\ref{thm:direct_vs_hierarchical} in Appendix~\ref{app:direct_vs_hierarchical}).

In our analysis of two-tier voting systems, we solve a previously held conjecture~\cite{boland1989majority,kao2019modular} that the lowest collective accuracy is reached if the number of voters per group is equal to the number of groups. Moreover, we generalize this finding to multi-tier hierarchical voting systems, proving that collective accuracy is minimized asymptotically when group sizes equal $N_{\rm d}^{1/n}$, where $N_{\rm d}$ is the total number of voters in a $n$ level sytem.

For heterogeneous voter competencies and abstention, we develop a generating function approach to describe such voting systems analytically without relying on approximations that have been used in earlier studies~\cite{hoeffding1956distribution,hodges1960poisson,boland1989majority}. Utilizing this framework, we provide an example illustrating that hierarchical voting systems are able to correct for the underrepresentation of voters with high competencies in groups or jurisdictions with high abstention rates.

We see a number of directions for future research. First, it would be interesting to investigate a dual outcome model in which voters receive utility from both the total outcome and their local outcome (that is, the within-group winner). Such a model could be used in both the standard Condorcet framework where there exists one ''correct'' outcome, or it could include an extension to allow for heterogeneous preferences over two competing outcomes. Second, future research may extend the model to account for correlated information in the proposed analytical accuracy analysis. Third, it would be fruitful to test the empirical implications by examining electoral volatility, the quality of representation, or voter satisfaction as a function of different indirect electoral designs.
\acknowledgements
LB acknowledges financial support from the SNF (P2EZP2\_191888) and thanks Wendelin Werner for his lecture on ``Randomness and Stability'' at GYSS 2020 that inspired parts of this work. The authors thank Malte Henkel, PJ Lamberson, Josh LeClair, and Scott E. Page for helpful discussions, as well as Xiaofeng Lin for research assistance.
\bibliography{refs.bib}
\bibliographystyle{naturemag}
\appendix
\section{Direct versus hierarchical voting systems}
\label{app:direct_vs_hierarchical}
To formulate our proof of Theorem~\ref{thm:direct_vs_hierarchical}, let $k=2 k'+1$ and $N_{\rm d} = 2 N_{\rm d}'+1$, such that
\begin{align}
\begin{split}
p_{\rm d}(\epsilon)\coloneqq p_{\rm d}'(\epsilon;N_{\rm d})  &=\sum\limits_{j=N_{\rm d}' + 1}^{2 N_{\rm d}'+1} {2 N_{\rm d}'+1 \choose j} \epsilon^j{(1-\epsilon)}^{2 N_{\rm d}'+1-j}\\
&=(1 - \epsilon)^{N_{\rm d}'} \epsilon^{N_{\rm d}'+1} {2 N_{\rm d}' + 1 \choose N_{\rm d}'+1}\, {}_2F_1\left(1,-N_{\rm d}';N_{\rm d}'+2;\frac{\epsilon}{\epsilon-1}\right)
\label{eq:recursive_1}
\end{split}
\end{align}
and
\begin{align}
\begin{split}
p_{i+1}(\epsilon) \coloneqq p_{\rm i+1}'(\epsilon;k') &=g(p_i(\epsilon))=\sum\limits_{j=k'+1 }^{2 k' + 1} {2 k'+1 \choose j} p^j_{i}{(1-p_{i})}^{2 k '+1-j},\\
&=(1 - p_i)^{k'} p_i^{k'+1} {2 k' + 1 \choose k'+1} \,{}_2F_1\left(1,-k';k'+2;\frac{p_i}{p_i-1}\right)\,,
\label{eq:recursive_2}
\end{split}
\end{align}
where we identified the truncated binomial sum with the function $g$, and ${}_2F_1$ is the (ordinary) hypergeometric function. Point (i) (i.e., $p_{\rm d}(\epsilon) = p_{n}(\epsilon)$ for  $\epsilon \in \{ 0,0.5,1\}$) can be readily obtained from the definitions of $p_{\rm d}(\epsilon)$ and $p_{i+1}(\epsilon)$. 
\begin{figure}
    \centering
    \includegraphics{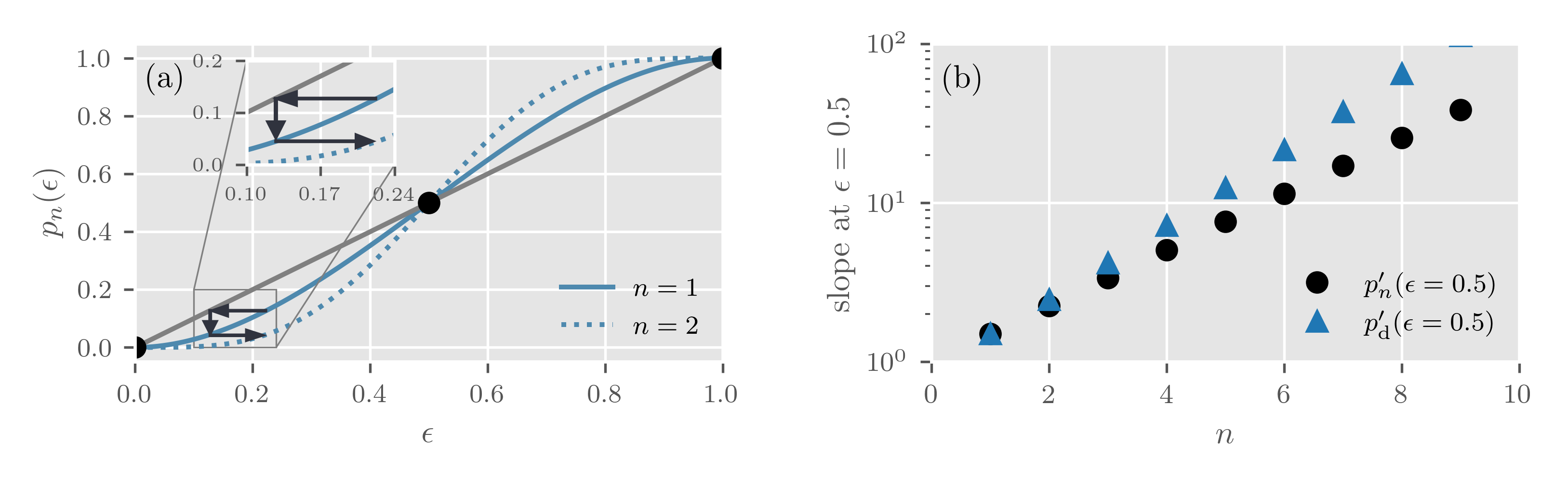}
    \caption{\textbf{Fixed-point iteration and hierarchical voting.} (a) The probability $p_n(\epsilon)$ of reaching the correct voting outcome in a hierarchical voting system as a function of $\epsilon$, the probability that voters in the electorate vote for the correct outcome. One can obtain $p_{n+1}(\epsilon)$ from $p_{n}(\epsilon)$ via an iteration (black arrows) from $n=1\rightarrow n+1=2$. Black disks indicate fixed points of the iteration. (b) The slope of $p_n(\epsilon)$ (black disks) and $p_{\rm d}(\epsilon)$ (blue triangles) at $\epsilon=0.5$. For $n\geq 2$, $p_{\rm d}'(\epsilon=0.5)$ is larger than $p_{n}'(\epsilon=0.5)$.}
    \label{fig:fixed_point}
\end{figure}
\begin{itemize}
\item For $\epsilon=0$, we find $p_{\rm d}(\epsilon)=p_n(\epsilon)=0$ (note that ${}_2F_1\left(a,b;c;z=0\right)=0$).
\item For $\epsilon=0.5$, we obtain $p_{1}(\epsilon)=(1-0.5)^{k'} 0.5^{k'+1} 4^{k'} = 0.5 $ and thus $p_i(\epsilon) = 0.5$ for all $i\in\{1,\dots,n\}$. Similarly, we find $p_{\rm d}(\epsilon) = 0.5$ for $\epsilon = 0.5$ and conclude that $p_{\rm d}(\epsilon) = p_n(\epsilon)=0.5$ for $\epsilon=0.5$.
\item For $\epsilon=1$, the last term in the sums \eqref{eq:recursive_1} and \eqref{eq:recursive_2} dominates, i.e., $p_{\rm d}(\epsilon) = p_n(\epsilon)=1$. Alternatively, we can use the identity
\begin{align}
\begin{split}
&\lim_{\epsilon \rightarrow 1} \, (1-\epsilon)^{k'} {}_2F_1\left(1,-k';k'+2;\frac{\epsilon}{\epsilon-1}\right)=
\lim_{\epsilon \rightarrow 1} \, (1-\epsilon)^{k'} (1-\epsilon)^{-k'} \left[ \frac{k'! \, (k'+1)!}{(2k'+1)!}+ \mathcal{O}(\epsilon-1)^2\right]\\
&=\frac{k'! \, (k'+1)!}{(2k'+1)!}={2 k' + 1 \choose k' + 1}^{-1}\,.
\end{split}
\end{align}
to show that $p_n(\epsilon)=1$ (and similarly $p_{\rm d}(\epsilon)=1$). 
\end{itemize}
A graphical interpretation of $p_{i+1}(\epsilon)=p_{i}(\epsilon)$ for $\epsilon\in\{0,0.5,1\}$ is that 0, 0.5, and 1 are fixed points of the iteration $p_{i+1}(\epsilon)=g(p_{i}(\epsilon))$ (black disks in Fig.~\ref{fig:fixed_point}a).

To prove points (ii--iii), we will use that $p_{\rm d}(\epsilon)$ and $p_{n}(\epsilon)$ are monotonically increasing with $\epsilon$ and convex for $\epsilon\in (0,0.5)$~\cite{boland1989majority}. For $p_{\rm d}(\epsilon)$ this follows directly from the definition \eqref{eq:recursive_1} and for $p_{n}(\epsilon)$ this follows from the way it is constructed via an iteration from $p_0(\epsilon)$ as we illustrate in Fig.~\ref{fig:fixed_point}(a). Everything that is left to show is that $p_{n}'(\epsilon)$ is smaller than $p_{\rm d}'(\epsilon)$ at $\epsilon=0.5$. As pointed out by Berg in \cite{berg1997indirect}, ``the derivatives at this point are crucial when two such functions are compared.'', motivating our following comparison of $p_{n}'(\epsilon)$ and $p_{\rm d}'(\epsilon)$. The derivative of $p_{\rm d}(\epsilon)$ with respect to $\epsilon$ at $\epsilon=0.5$ is
\begin{align}
\begin{split}
\partial_{\epsilon} p_{\rm d}(\epsilon=0.5)&=0.25^{N_{\rm d}'-1} \frac{(2 N_{\rm d}' + 1)!}{N_{\rm d}'!}\left[0.25\, \frac{{}_2F_1\left(1,-N_{\rm d}';N_{\rm d}'+2;-1\right)}{(N_{\rm d}'+1)!}+0.5\, N_{\rm d}'\,\frac{{}_2F_1\left(2,1-N_{\rm d}';N_{\rm d}'+3;-1\right)}{(N_{\rm d}'+2)!} \right]\\
&=0.25^{0.5(k^n-1)} \frac{(k^n)!}{(0.5(k^n-1))!}\left[ \frac{{}_2F_1\left(1,0.5(1-k^n);0.5(k^n+3);-1\right)}{(0.5(k^n+1))!}\right.\\
&\left.+(k^n-1)\,\frac{{}_2F_1\left(2,0.5(3-k^n);0.5(k^n+5);-1\right)}{(0.5(k^n+3))!} \right]\,.
\end{split}
\label{eq:pdp_05}
\end{align}
For $p_0(\epsilon)=g(\epsilon)$, we find
\begin{align}
\begin{split}
\partial_{\epsilon} g(\epsilon=0.5)&=0.25^{k'-1} \frac{(2 k' + 1)!}{k'!}\left[0.25\, \frac{{}_2F_1\left(1,-k';k'+2;-1\right)}{(k'+1)!}+0.5\, k'\,\frac{{}_2F_1\left(2,1-k';k'+3;-1\right)}{(k'+2)!} \right]\\
&=0.25^{0.5(k-1)} \frac{k!}{(0.5(k-1))!}\left[ \frac{{}_2F_1\left(1,0.5(1-k);0.5(k+3);-1\right)}{(0.5(k+1))!}\right.\\
&+\left.(k-1)\,\frac{{}_2F_1\left(2,0.5(3-k);0.5(k+5);-1\right)}{(0.5(k+3))!} \right]\,.
\end{split}
\label{eq:derivative_g}
\end{align}
The derivative of $p_{n}(\epsilon)$ is
\begin{equation}
\partial_{\epsilon} p_{n}(\epsilon)= \partial_{\epsilon}(\underbrace{g\circ g \ldots g \circ g}_\text{$n$ times}) \circ \epsilon=\underbrace{(\partial_{\epsilon}g \circ \epsilon)\, (\partial_{\epsilon} g\circ g \circ \epsilon)\ldots (\partial_{\epsilon} g\circ \underbrace{g \circ \ldots \circ g}_{\text{$n-1$ times}} \circ \epsilon)}_{\text{$n$ times}}\,.
\end{equation}
Using $g(\epsilon=0.5)=0.5$ yields
\begin{align}
\begin{split}
\partial_{\epsilon}p_{n}(\epsilon=0.5) &= \left[\partial_{\epsilon} g(\epsilon=0.5)\right]^n\\
&=0.25^{0.5n(k-1)} \frac{(k!)^n}{\left[(0.5(k-1))!\right]^n}\left[ \frac{{}_2F_1\left(1,0.5(1-k);0.5(k+3);-1\right)}{(0.5(k+1))!}\right.\\
&+\left.(k-1)\,\frac{{}_2F_1\left(2,0.5(3-k);0.5(k+5);-1\right)}{(0.5(k+3))!} \right]^n\,.
\end{split}
\label{eq:pnp_05}
\end{align}
Note that $\partial_{\epsilon} p_{\rm d}(\epsilon=0.5)=\partial_{\epsilon} p_{n}(\epsilon=0.5)$ for $n=1$ and arbitrary $k$. According to Eq.~\eqref{eq:pnp_05}, $\log \left[\partial_{\epsilon} p_{n}(\epsilon=0.5)\right]$ grows linearly with $n$ (with slope $\log\left[\partial_{\epsilon} g(\epsilon=0.5)\right]$). Because of the $(k^n)!$ term in Eq.~\eqref{eq:recursive_1}, the slope of $\partial_{\epsilon} p_{\rm d}(\epsilon=0.5)$ for constant $k$ is larger than that of $\partial_{\epsilon} p_{n}(\epsilon=0.5)$. Instead of the linear growth of $\log \left[\partial_{\epsilon} p_{n}(\epsilon=0.5)\right]$ with $n$, we find that $\log\left[\partial_{\epsilon} p_{\rm d}(\epsilon=0.5)\right]$ involves a $\log\left[(k^n)!\right]= n k^n \log(k)-k^n + \mathcal{O}\left[\log(k^n)\right]$ term. In Fig.~\ref{fig:fixed_point}(b), we show $\partial_{\epsilon} p_{\rm d}(\epsilon=0.5)$ and $\partial_{\epsilon} p_{n}(\epsilon=0.5)$ for $k=3,5$ as a function of $n$, confirming that $\partial_{\epsilon} p_{\rm d}(\epsilon=0.5)$ is larger than $\partial_{\epsilon} p_{n}(\epsilon=0.5)$ for $n\geq 2$. Given the discussed properties of $p_{\rm d}(\epsilon)$ and $p_n(\epsilon)$, we have thus shown that for $n\geq 2$, $p_{\rm d}(\epsilon) < p_{n}(\epsilon)$ for $\epsilon \in (0,0.5)$ (point ii) and $p_{\rm d}(\epsilon) > p_{n}(\epsilon)$ for $\epsilon \in (0.5,1.0)$ (point iii).
\section{Group-size associated with lowest performance}
\label{app:lowest_accuracy}
\begin{figure}
\includegraphics{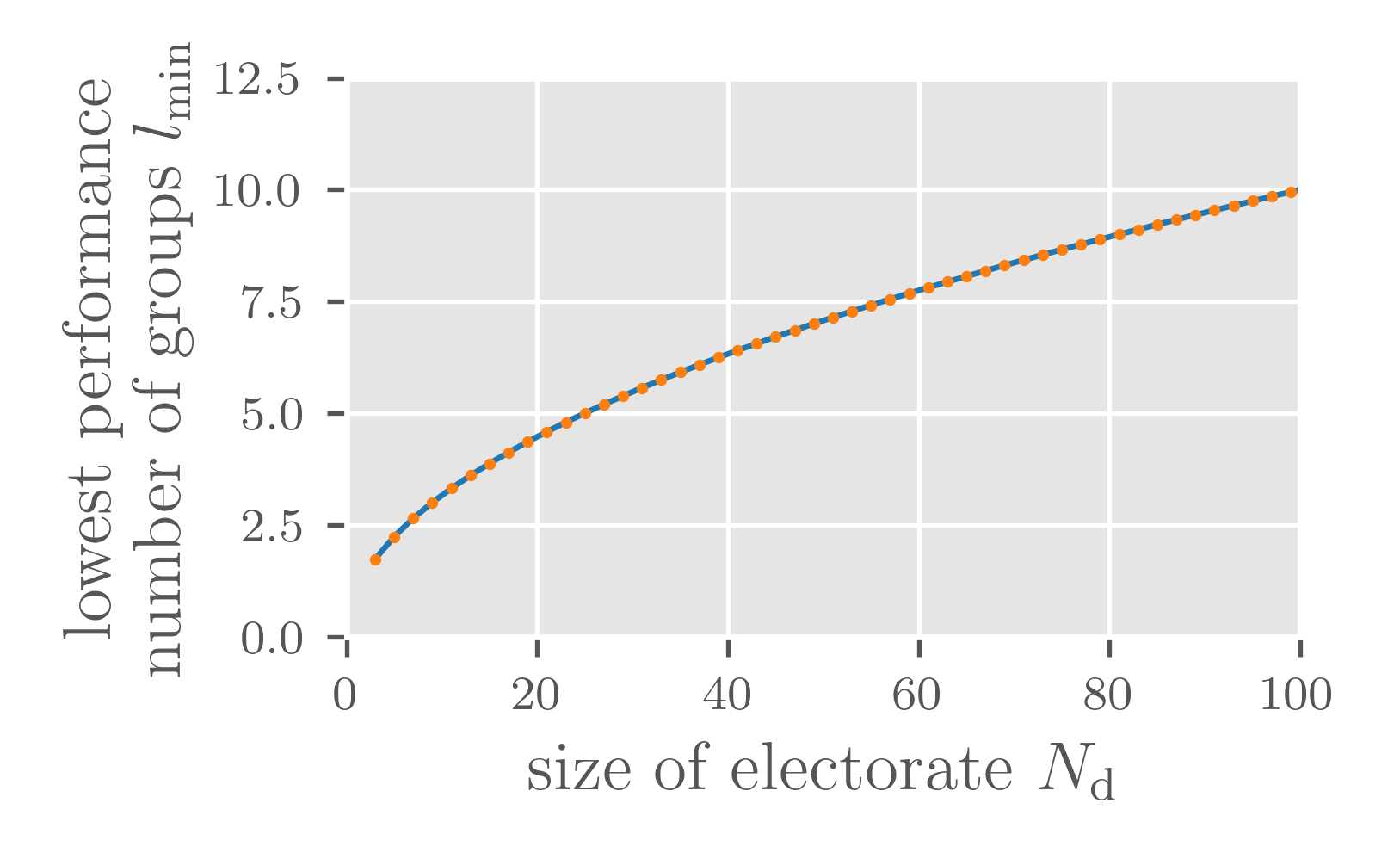}
    \caption{\textbf{Number of groups with the lowest performance as a function of the electorate size.} Orange dots were obtained by numerically determining the minimum of the product \eqref{eq:derivative_exact}. Note that Eq.~\eqref{eq:derivative_exact} admits real-valued solutions of $l_{\rm min}$. The solid blue line describes the number of groups with the lowest performance according to $l_{\rm min}=\sqrt{N_{\rm d}}$.}
\label{fig:minimum_slope}
\end{figure}
We aim at determining the values of $k$ (number of voters per group), $l$ (number of groups) that are associated with the lowest collective accuracy in a two-hierarchy system [i.e., with the smallest value of $p_2(\epsilon;k)$ for $\epsilon>0.5$ in Eq.~\eqref{eq:p2_2tier}]. As in Appendix \ref{app:direct_vs_hierarchical}, we use $k'=(k-1)/2$ and $l'=(l-1)/2$. Note that $l'$ is a function of $k'$ for fixed $N=kl$.

To determine the values of $k'$ and $l'$, and hence of $k$ and $l$, that are associated with the smallest values of $p_2(\epsilon;k)$ for $\epsilon>0.5$, we study the slope of $p_2(\epsilon;l')$ as a function of $l'$ in the vicinity of $\epsilon=0.5$. Since the considered reliability functions are concave increasing in the interval $[1/2,1]$~\cite{boland1989majority,berg1997indirect}, it is sufficient to evaluate of $p_2(\epsilon;l')$ for a fixed value of $\epsilon$ and varying $l'$. An expansion of $p_2(\epsilon;l')$ about $(\epsilon,l')=(0.5+\Delta \epsilon,l')$ with $\Delta \epsilon>0$ yields
\begin{align}
\begin{split}
&p_2(\epsilon=0.5+\Delta \epsilon;l')=0.5+\Delta \epsilon \left.\frac{\partial p_2}{\partial \epsilon}\right|_{(\epsilon,l')=(0.5,l')}+\mathcal{O}({(\Delta \epsilon)}^2)\,.
\end{split}
\end{align}
The derivative $\partial_{\epsilon}p_2$ (or Banzhaf number~\cite{berg1997indirect}) is
\begin{equation}
\frac{\partial p_2}{\partial \epsilon}(p_1(\epsilon=0.5;k');l')=\frac{\partial p_1}{\partial \epsilon}(\epsilon=0.5;k')\,\frac{\partial p_2}{\partial \epsilon}(p_1(\epsilon=0.5;k');l')=\frac{\partial p_1}{\partial \epsilon}(\epsilon=0.5;k')\,\frac{\partial p_2}{\partial \epsilon}(\epsilon=0.5;l')\,.
\end{equation}
As a side note, observe that it is symmetric in $(k,l)=(k,N_{\rm d}/k)$. This symmetry in $(k,l)$ is a consequence of the fixed-point behavior of $p_2$ at $\epsilon=0.5$ (black disks in Fig.~\ref{fig:fixed_point}a). 

For evaluating $\partial_{\epsilon} p_2$, we use Eq.~\eqref{eq:derivative_g}, which we rewrite to obtain
\begin{align}
\begin{split}
\partial_{\epsilon} g(\epsilon=0.5;k')=1+\frac{4k' \Gamma(k'+3/2)\,{}_2F_1\left(2,1-k';k'+3;-1\right)}{\sqrt{\pi}\Gamma(k'+3)} \,,
\end{split}
\label{eq:derivative_g2}
\end{align}
and find
\begin{equation}
\frac{\partial p_2}{\partial \epsilon}(p_1(\epsilon=0.5;k');l')=\partial_{\epsilon} g(\epsilon=0.5;k') \partial_{\epsilon} g(\epsilon=0.5;l')\,.
\label{eq:derivative_exact}
\end{equation}
Invoking Stirling's approximation in the limit of large $k'$, we approximate the ratio of the Gamma functions in Eq.~\eqref{eq:derivative_g2} according to
\begin{equation}
\frac{\Gamma(k'+3/2)}{\Gamma(k'+3)} \sim {k'}^{-\frac{3}{2}}\,.
\end{equation}
To derive the corresponding asymptotic relation for the hypergeometric function ${}_2F_1\left(2,1-k';k'+3;-1\right)$ in Eq.~\eqref{eq:derivative_g2}, we use relation 7.3.6.3 of \cite{prudnikov} and obtain
\begin{equation}
\lim_{k'\rightarrow\infty}{}_2F_1\left(2,1-k';k'+3;-1\right)= \frac{k'}{2}-\frac{\sqrt{\pi}}{4}\sqrt{k'}+\frac{3}{2}+\mathcal{O}(k'^{-\frac{1}{2}})\,.
\end{equation}
We thus find
\begin{equation}
\partial_{\epsilon} g(\epsilon=0.5;k')\sim \frac{2 (k'+3)}{\sqrt{\pi k'}}\,.  
\label{eq:derivative_asymptotic}
\end{equation}
For large $k$ and $l$, we set $k'\sim k/2$, $l'\sim N_{\rm d}/(2k)$ and obtain
\begin{equation}
\frac{\partial p_2}{\partial \epsilon}(p_1(\epsilon=0.5;k/2);N_{\rm d}/(2k))\sim\frac{2 (k+6) (6 k+N_{\rm d})}{\pi  k \sqrt{N_{\rm d}}}\,.
\label{eq:asymptotic_derivative}
\end{equation}
The minimum of Eq.~\eqref{eq:asymptotic_derivative} is attained for 
\begin{equation}
(k_{\rm min},l_{\rm min})=(\sqrt{N_{\rm d}},\sqrt{N_{\rm d}})\,.
\label{eq:minimum}
\end{equation}
We thus conclude that the minimum of the first derivative of $p_2(\epsilon;l)$ in the vicinity of $\epsilon=0.5$ is attained asymptotically for $l_{\rm min}=k_{\rm min}=\sqrt{N_{\rm d}}$. Hence, $p_2(\epsilon;l=\sqrt{N_{\rm d}})$ corresponds to the reliability function with the smallest collective accuracy for $\epsilon>0.5$ and the largest collective accuracy for $\epsilon<0.5$.

We derived Eq.~\eqref{eq:minimum} using different asymptotic relations which are valid for large $k'$ and $l'$. Based on numerical calculations, we observe that the minimum of the product \eqref{eq:derivative_exact} follows the same square-root law for small values of $k'$ and $l'$ (Fig.~\ref{fig:minimum_slope}).

The above proof can be extended to general hierarchical systems with groups of sizes $k^{(1)},k^{(2)},\dots,k^{(n)}$ in layers $1,2,\dots,n$, respectively. The total number of voters in the electorate is $N_{\rm d}=k^{(1)}k^{(n)}\cdots k^{(n)}$. In accordance with Eqs.~\eqref{eq:derivative_exact} and \eqref{eq:derivative_asymptotic}, observe that the reliability function $p_n(\epsilon)$ of a hierarchical voting system with $n$ layers factorizes at the fixed point $\epsilon=0.5$ and is given by
\begin{align}
\partial_{\epsilon} p_n(\epsilon=0.5) &= \prod_{i=1}^n \partial_{\epsilon} g(\epsilon=0.5;k'^{(i)}) \\
&\sim \prod_{i=1}^n \frac{\sqrt{2}(k^{(i)}+6)}{\sqrt{\pi k^{(i)}}}\,,
\label{eq:multi-tier-performance}
\end{align}
where $k'^{(i)}\sim k^{(i)}/2$ and $k^{(n)}=N_{\rm d}/(k^{(1)} \cdots k^{(n-1)})$. Determining the minimum of Eq.~\eqref{eq:multi-tier-performance} yields the following $n-1$ equations with $n-1$ unknowns:
\begin{align}
\begin{split}
k_{\rm min}^{(1)} \, k_{\rm min}^{(2)} \, \cdots \, k_{\rm min}^{(n-2)}\,  {k_{\rm min}^{(n-1)}}^2 &= N_{\rm d}\\
k_{\rm min}^{(1)} \, k_{\rm min}^{(2)} \, \cdots\,  {k_{\rm min}^{(n-2)}}^2 \, k_{\rm min}^{(n-1)} &= N_{\rm d}\\
&\vdots\\
k_{\rm min}^{(1)} \, {k_{\rm min}^{(2)}}^2 \, \cdots \, k_{\rm min}^{(n-2)} \, k_{\rm min}^{(n-1)} &= N_{\rm d}\\
{k_{\rm min}^{(1)}}^2 \,  k_{\rm min}^{(2)}\,  \cdots \, k_{\rm min}^{(n-2)}\,  k_{\rm min}^{(n-1)} &= N_{\rm d}\,.
\end{split}
\label{eq:general_minimum}
\end{align}
The solution of the above set of equations is $k_{\rm min}^{(1)}=k_{\rm min}^{(2)}=\cdots=k_{\rm min}^{(n)}=N_{\rm d}^{1/n}$.

\section{The fewest voters needed to sway an outcome occurs at $\sqrt{N_{\rm d}}$}
\label{app:fewest_votrs}
In a two-tier hierarchical system, the lowest number of voters needed to support a winning measure is equal to 
\begin{equation}
\left(\frac{k+1}{2}\right)\left(\frac{\frac{N_{\rm d}}{k}+1}{2}\right)=\frac{N_{\rm d}+k+\frac{N_{\rm d}}{k}+1}{4}\,,
\label{eq:2tier_voters}
\end{equation}
where $k$ denotes the odd number of voters in each of the $N_{\rm d}/k$ groups.

Taking the derivative of Eq.~\eqref{eq:2tier_voters} with respect to $k$ gives
\begin{equation}
    \frac{1}{4}-\frac{N_{\rm d}}{4k^2}\,.
\end{equation}
Setting this to zero and solving for $k$ reveals that $k_{\rm min}=\sqrt{N_{\rm d}}$. Given the fact that $k_{\rm min}$ must be a positive number, this implies a unique extremum. Taking the second derivative of Eq.~\eqref{eq:2tier_voters} yields
\begin{equation}
    \frac{N_{\rm d}}{2k^3}\,.
\end{equation}
This number is positive for $k_{\rm min}=\sqrt{N_{\rm d}}$; therefore, it is a unique minimum.\\

More generally, in an $n$ level hierarchy with $k^{(i)}$ groups (voters) at each level $i=1,\ldots,n$, where $k^{(i)}$ is odd for all $i$, the fewest voters needed to support a winning measure can be found by searching for the minimum of
\begin{equation*}
  f(k^{(1)},\ldots,k^{(n-1)})=
  \left(\frac{k^{(1)}+1}{2}\right)\left(\frac{k^{(2)}+1}{2}\right)\cdots\left(\frac{\frac{N_{\rm d}}{k^{(1)} k^{(2)}\cdots k^{(n-1)}}}{2}\right).
\end{equation*}
For each $i\in\{1,\ldots,n-1\}$ the partial derivative with respect to $k^{(i)}$ is
\begin{equation*}
  \frac{\partial }{\partial k^{(i)}}f(k^{(1)},\ldots,k^{(n-1)})=\frac{\left(\prod_{j\neq i}(k^{(j)}+1)\right)\left(k^{(i)}\left(\prod_{j=1}^{n-1}k^{(j)}\right)-N_{\rm d}\right)}{2^n k^{(i)}\prod_{j=1}^{n-1}k^{(j)}}.
\end{equation*}
The unique positive solution to $\nabla f(k_{\rm min}^{(1)},\ldots,k_{\rm min}^{(n-1)})=(0,\ldots,0)$ is $k_{\rm min}^{(i)}=\sqrt[n]{N_{\rm d}}$ and can be found by setting $k_{\rm min}^{(i)}\left(\prod_{j=1}^{n-1}k_{\rm min}^{(j)}\right)=N_{\rm d}$. Observe that the resulting set of equations is equivalent to Eq.~\eqref{eq:general_minimum}. The eigenvalues of the Hessian at this solution are
\begin{equation*}
  \left\{\frac{n\left(1+\sqrt[n]{N_{\rm d}}\right)^{n-2}}{2^n \sqrt[n]{N_{\rm d}}},\frac{\left(1+\sqrt[n]{N_{\rm d}}\right)^{n-2}}{2^n \
  \sqrt[n]{N_{\rm d}}}\right\}.
\end{equation*}
Since these are both positive, the Hessian is positive definite, so $k^{(i)}_{\rm min}=\sqrt[n]{N_{\rm d}}$ for $i=1,\ldots,n$ is a minimum of $f$.
\section{Generating function approach}
\label{sec:generating_function}
Let $X_i$ be a Bernoulli random variable for which $\Pr(X_i=1)=p_i$ and $\Pr(X_i=0)=q_i=1-p_i$. In what follows, we assume that $p_i\neq 1$. To derive the probability mass function of the number of positive outcomes
\begin{equation}
\Sigma_n\equiv \sum_{i=1}^n X_i\,,
\end{equation}
we start from the generating function~\cite{percus1985probability}
\begin{equation}
G_n(t)=\sum_{l=0}^\infty t^l P(\Sigma_n=l)\,,    
\end{equation}
which converges for $|t|\leq 1$. Using that $G_n(t)=\mathbb{E}(t^{\Sigma_n})=\mathbb{E}(t^{\sum_i X_i})=\prod_i \mathbb{E}(t^{X_i})=\prod_i \left(q_i+p_it\right)$, the generating function can be written as~\cite{percus1985probability}
\begin{equation}
G_n(t)=\prod_{i=1}^n q_i \prod_{i=1}^n \left(1+\frac{p_i}{q_i} t\right)\,.
\label{eq:generating_function}
\end{equation}
The probability that $\Sigma_n$ is equal to $s\leq n$ is
\begin{align}
\begin{split}
\Pr(\Sigma_n=s)&=\frac{1}{s\,!} \left.\frac{\partial^s G_n(t)}{\partial t^s}\right|_{t=0}\\
&=\prod_{i=1}^n q_i \frac{1}{s!}\sum_{k_1+k_2+\dots+k_n=s} {s \choose k_1, k_2, \dots,k_n }\prod_{1\leq i \leq n} \left.\frac{\partial^{k_i}}{\partial t^{k_i}}\left(1+\frac{p_i}{q_i} t\right)\right|_{t=0}\\
&=\prod_{i=1}^n q_i \sum_{1\leq k_1 < k_2 <\dots<k_n\leq s} \frac{p_{k_1}}{q_{k_1}}\dots\frac{p_{k_n}}{q_{k_n}}\\
&=\frac{1}{P} C_s(p_1/q_1,\dots,p_n/q_n)\,,
\end{split}
\label{eq:pmf}
\end{align}
where $P=1/\prod_i q_i$ and $C_s$ is the $s$th elementary symmetric polynomial. In the second step of Eq.~\eqref{eq:pmf}, we used the general Leibniz rule to express the $s$th derivative of the second product in Eq.~\eqref{eq:generating_function} with respect to $t$.

To calculate the $s$th elementary symmetric polynomial, we utilize the recursion relation (``Newton's identities'')~\cite{mead1992newton}
\begin{align}
C_s&=\frac{1}{s} \sum_{k=0}^{s-1}(-1)^k C_{s-1-k}\sigma_{k+1}\,,\quad s>0\,,\\
C_0&=1\,,
\end{align}
where
\begin{equation}
\sigma_s(p_1/q_1,\dots,p_n/q_n)=\sum_{i=1}^n \left(\frac{p_i}{q_i}\right)^s\,.
\end{equation}

To apply Eq.~\eqref{eq:pmf} to the voting system with abstention and heterogeneous opinions (see main text), positive outcomes ($X_i=1$) have to be in the majority. For $n=2 n'+1$, positive outcomes are in the majority if $\Sigma_n \geq n'+1$. The corresponding probability is
\begin{equation}
\Pr(\Sigma_n \geq n'+1)=\sum_{s=n'+1}^n\Pr(\Sigma_n = s)\,.
\end{equation}
As an example, we consider the case with $n=3$ independent Bernoulli random variables. Positive outcomes are in the majority if $\Sigma_3\geq 2$. Using the elementary symmetric polynomials for $n=2$ and $s=2,3$ (Tab.~\ref{tab:Csn}), the probability of $\Sigma_3\geq 2$ is
\begin{equation}
\Pr(\Sigma_3 \geq 2)=\frac{1}{P}\left(\frac{p_{1} p_{2}}{q_{1} q_{2}}+\frac{p_{1} p_{3}}{q_{1} q_{3}}+\frac{p_{2} p_{3}}{q_{2} q_{3}}+\frac{p_{1} p_{2} p_{3}}{q_{1}q_{2}q_{3}}\right)\,,
\end{equation}
where $P=1/(q_1 q_2 q_3)$.
\begin{table}[htp!]
\centering
\renewcommand*{\arraystretch}{2}
\begin{tabular}{|>{\centering\arraybackslash} m{4em}|>{\centering\arraybackslash} m{8em}|>{\centering\arraybackslash} m{10.5em}|>{\centering\arraybackslash} m{5em}|}
\hline
\backslashbox{$\,\,\,n\,\,\,$}{$\,\,\,s\,\,\,$} & $1$ & $2$ & $3$ \\ \hline\hline
$1$ & $\displaystyle{p_1\over q_1}$ & - & - \\\hline
$2$ & $\displaystyle{{p_1\over q_1}+{p_2\over q_2}}$ & $\displaystyle{p_1 p_2\over q_1 q_2}$ & - \\\hline
$3$ & $\displaystyle{{p_1\over q_1}+{p_2\over q_2}+{p_3\over q_3}}$ & $\displaystyle{{p_1 p_2\over q_1 q_2}+{p_1 p_3\over q_1 q_3}+{p_2 p_3\over q_2 q_3}}$ & $\displaystyle{{p_1 p_2 p_3 \over q_1 q_2 q_3}}$\\\hline
\end{tabular}
\caption{Elementary symmetric polynomials $C_s(p_1/q_1,\dots,p_n/q_n)$ for $1\leq n \leq 3$ and $1\leq s\leq n$.}
\label{tab:Csn}
\end{table}
\end{document}